\documentclass[journal]{IEEEtran}

\usepackage{amsmath,epsfig,amssymb,amsfonts, psfrag}
\usepackage{graphicx}
\usepackage{color}
\usepackage{times}
\usepackage[lined,boxed,commentsnumbered]{algorithm2e}
\usepackage[mathscr]{eucal} 
\usepackage{tikz}
\usepackage{pgfplots}
\usepackage{epstopdf} 

\usetikzlibrary{arrows}

\newcounter{theorem}
\newtheorem{theorem}{Theorem}[section]
\newtheorem{lemma}[theorem]{Lemma}
\newtheorem{corollary}[theorem]{Corollary}
\newtheorem{proposition}[theorem]{Proposition}

\newtheorem{conjecture}[theorem]{Conjecture}

\newenvironment{proof}{ \vspace{\smallskipamount}\par {\it Proof.} }{\hspace*{\fill}~\IEEEQED\par} 

\newcommand{\set}[1]{\mathscr{#1}} 

\newcommand{\tens}[1]{\boldsymbol{\mathcal{#1}}}

\newcommand{\matr}[1]{\boldsymbol{#1}}

\newcommand{\vect}[1]{\boldsymbol{#1}}

\newcommand{\T}{{\sf T}}      
\newcommand{\krank}[1]{\mathop{\operator@font krank}\{#1\}}
\newcommand{\trace}[1]{\mathop{\operator@font trace}\{#1\}}
\newcommand{\Diag}[1]{\mathop{\operator@font Diag}\{#1\}}    
\newcommand{\diag}[1]{\mathop{\operator@font diag}\{#1\}}    
\newcommand{\Span}[1]{\mathop{\operator@font Span}\{#1\}}    
\newcommand{\out}{\otimes}    
\newcommand{\khatri}{\odot}   
\newcommand{\hadam}{\boxdot}  
\newcommand{\kron}{\boxtimes} 
\makeatletter
\newcommand{\rank}[1]{\mathop{\operator@font rank}\{#1\}}   
\renewcommand{\vec}{\mathop{\operator@font vec}}   %
\newcommand{\Unvec}{\mathop{\operator@font Unvec}}   %
\makeatother

\newcommand{\KK}{\mbox{$\mathbb{K}$}} 

\definecolor{brickred}{rgb}{0.6,0,0}

\definecolor{brickblue}{rgb}{0,0,0.6}

\definecolor{brickgreen}{rgb}{0,0.5,0}

\begin{document}\sloppy

\title{ Rank-1 Tensor Approximation Methods\\ and Application to Deflation }

\author{Alex~P.~da~Silva,~
       Pierre~Comon,~\IEEEmembership{Fellow,~IEEE,}
        and~Andr\'e~L.~F.~de~Almeida,~\IEEEmembership{Senior Member,~IEEE}
\thanks{This work has been funded by the European Research Council under the 7th Framework  Program  FP7/2007-2013 Grant Agreement no. 320594, and  by the Conselho Nacional de Desenvolvimento Cient\'ifico e Tecnol\'ogico (CNPQ) under the program Ci\^encias sem Fronteiras.}}

\markboth{Journal of \LaTeX\ Class Files}%
{Shell \MakeLowercase{\textit{et al.}}: Bare Demo of IEEEtran.cls for Journals}

\maketitle

\begin{abstract}
Because of the attractiveness of the canonical polyadic (CP) tensor decomposition in various applications,  several algorithms have been designed to compute it, but  efficient ones are still lacking. Iterative deflation algorithms based on successive rank-1 approximations can be used to perform this task, since the latter are rather easy to compute.  We first present an algebraic rank-1 approximation method that performs better than the standard  higher-order singular value decomposition (HOSVD)  for three-way tensors. Second, we propose a new iterative rank-1 approximation algorithm that improves any other rank-1 approximation method. Third, we describe a probabilistic framework allowing to  study  the convergence of deflation CP decomposition (DCPD) algorithms based on successive rank-1 approximations. A set of computer experiments then validates theoretical  results  and demonstrates the efficiency of DCPD algorithms compared to other ones.
\end{abstract}


\section{Introduction}\label{sec:intro}

 In the last years, tensors have been playing an important role in many applications such as blind source separation \cite{ComoJ10, castella2007blind}, telecommunications \cite{AlmeFM07:SP}, chemometrics \cite{smilde2005multi}, neuroscience \cite{becker2014eeg}, sensor array processing \cite{sahnounjoint} and data mining \cite{Sava08:thesis}. The attractiveness behind tensors lies in the uniqueness of their canonical polyadic (CP) decomposition  under mild conditions \cite{Kier00:jchemo}, which is a powerful property not shared by standard matrix-based tools. There are several methods to compute the CP tensor decomposition.  We  will point out here some of the most used methods among many others. For the exact CP decomposition, \cite{ten1991kruskal} proposes a direct computation to decompose $2\times n \times n$ tensors. In \cite{brachat2010symmetric}, a generalization of Sylvester's algorithm is described for decomposing symmetric tensors. In \cite{de2006link}, one can use simultaneous matrix diagonalization by congruence, provided that the rank of the tensor is smaller than its greatest dimension. An approach based on eigenvectors of tensors is proposed in \cite{oeding2013eigenvectors}.

In practice, tensors are corrupted by noise so that one needs to compute an approximate decomposition of given rank. 
 Computing the exact CP decomposition is difficult \cite{hillar2013most}, but  finding a lower-rank approximation is even harder. 
 In fact, this is an ill-posed problem in general \cite{de2008tensor}. Nevertheless, some useful algorithms have been conceived to solve locally the low-rank approximation problem. This kind of algorithms can be found in \cite{comon2009tensor, kolda2009tensor, rajih2008enhanced, sorber2013optimization, tomasi2006comparison}. One of the most widely used is the alternating least squares (ALS) algorithm \cite{smilde2005multi}, which is an iterative method that consists in conditionally updating in an alternate way, the matrix factors composing the CP decomposition.  Other gradient and Newton-based methods estimate the factor matrices all-at-once. Howsoever, all theses algorithms  have disappointing convergence properties \cite{comon2009tensor, Paat97:cils}. Another kind of algorithms is based on rank-1 deflation. It is known that the conventional deflation works for matrices but not for tensors \cite{StegC10:laa}. In  \cite{phantensor},   the authors propose deflation methods that work only if the rank of the tensor is not larger than its dimensions. In \cite{CichP09:ieice}, ALS is used to update the columns of matrix factors in a deflation procedure, but for non-negative tensors. However, these deflation methods strongly depend on initialization and no convergence study has been conducted. 

In the same vein of iterative deflation, the authors proposed in \cite{AlexIcassp2015}  a deflation-based CP decomposition (DCPD), based on successive rank-1 approximations computed by low-complexity methods.
 Other rank-1 approximation methods can be used in DCPD, for instance, ALS. However, the latter exhibits an unbounded complexity and no satisfactory convergence study is available for rank-1 approximation apart \cite{WangC14:simax}, which shows results on the global convergence for generic tensors in the sense that, for any initialization, ALS converges to the same point in general.  A quasi-Newton method defined on Grassmannians are developed in \cite{savas2010quasi}. However, they exhibit an unbounded complexity to compute rank-1 approximations, since these methods need to be iterated. 
 In \cite{lasserre2001global}, the best rank-1 approximation problem can be computed by means of an algebraic geometry moment method, but it is only applicable to very small tensor dimensions since the number of variables grows exponentially when building convex relaxations. Moreover, even for small dimensions, its convergence  is  very slow.  In \cite {bucero2014border},  the authors propose an improvement method for \cite{lasserre2001global} based on border basis relaxation, but again the method is limited to small tensor dimensions.  Semidefinite relaxations are proposed in \cite{nie2014semidefinite} to compute rank-1 approximations, however the  convergence  becomes very slow when  dimensions are large.

\smallskip

In this paper, we report mainly three contributions. First, we  propose  an algebraic rank-1 approximation method, namely the sequential rank-1 approximation and projection (SeROAP), which can perform better than the standard truncated high-order singular value decomposition (THOSVD) \cite{de2000multilinear}. Indeed, we prove that the rank-1 approximation performed by SeROAP is always better than  the one  obtained from THSOVD for three-way tensors. Moreover, for  large  dimensions and small orders, we show that the computational complexity of SeROAP  is dramatically lower than that of THOSVD. 

Second, we propose an  alternating eigenvalue rank-1 iterative algorithm for three-way tensors, namely CE (coupled-eigenvalue), that improves other rank-1 approximation algorithms. We prove that if the solution obtained from some rank-1 approximation algorithm (e.g., SeROAP, THOSVD) is the input of the CE algorithm, the performed rank-1 approximation remains the same in the worst case. We also prove that the convergence to a stationary point is always guaranteed. Actually, results have shown that when the initialization  of the CE algorithm  is close enough to the global solution,  it recovers the best rank-1 approximation. Furthermore, when one dimension is much larger than the other two dimensions, the computational complexity of the CE algorithm can be lower than that of the standard ALS algorithm.

Third, we perform a theoretical study on deflation in order to analyze the convergence of the DCPD algorithm. In a first stage,  we show that the norm of residuals is monotonically reduced within the iterative deflation process. We also prove that the DCPD algorithm recovers the exact CP decomposition of a given tensor when residuals do not fall within a cone with an arbitrary small volume. In a second  stage,  we prove that the iterative deflation method can reduce the norm of the initial residual by a factor smaller than $\left(\sin(\beta)\right)^{L-1}$ ($\beta$  being  the angle of a suitable cone where the residuals can fall in) after $L$ iterations with high probability, when tensors are distributed according to an absolutely continuous probability measure,  and the probability function of residuals is continuous on some suitable angular interval. We also present a conjecture stating the existence of probability measures ensuring the convergence of the DCPD algorithm to an exact decomposition with high probability.

The paper is organized as follows. In Section \ref{sec: dcpt_alg} some standard iterative methods and  the DCPD algorithm as well as SeROAP and THOSVD methods   are described.  The computational complexity per iteration for each  algorithm  is provided.
The next two sections form the core of the paper. In Section \ref{sec:rank_one_res}, we first  prove that SeROAP performs better than THOSVD as far as rank-1 approximations of three-way tensors is concerned.
Then in a second part,  the CE algorithm is presented as well as the proof that it can refine any other rank-1 approximation method  and the proof of its convergence. In Section \ref{sec:pfDCPD},  among other theoretical results, we study conditions ensuring the convergence of the DCPD algorithm.  Finally, in Section \ref{sec:sim_res}, computer results show satisfactory performances of the proposed DCPD and rank-1 approximation methods,  compared to other related algorithms,  even under noisy scenarios. 

\section{Notation}\label{sec:not}

The notation employed in this paper is the following. Scalar numbers are denoted by lowercase letters and vectors by boldface lowercase ones. For matrices and tensors, boldface capital and calligraphic letters are used, respectively. Plain capitals are used to denote array dimensions. Depending on the context, greek letters can denote either  scalars, vectors, matrices or tensors. The symbols $\hadam, \khatri, \kron$ and $\otimes$ denote the Hadamard, Khatri-Rao, Kronecker and tensorial products, respectively, and $^{+}$ denotes matrix pseudo-inversion.  
The Euclidean scalar product between tensors is denoted by $\langle \tens{T}, \tens{U} \rangle = \sum_{i_{1}\cdots i_{N}} T_{i_{1}\cdots i_{N}}U^*_{i_{1}\cdots i_{N}}$. 
The angle between two tensors $\tens{U}$ and $\tens{V}$ will  refer to $\arccos\{|\langle\tens{U},\tens{V}\rangle|\}/ \|\tens{U}\|_F\, \|\tens{V}\|_F \in [0,\,\pi/2]$. 
$\|\cdot\|$ then denotes the Frobenius norm induced by the previous scalar product. We shall also use the $\ell^2$ operator norm $\|\cdot\|_2$ for matrices, which corresponds to the largest singular value. 
The mode-$n$ unfolding of a tensor $\tens{T}$ is denoted as $\matr{T}_{(n)}$, as proposed in \cite{kolda2009tensor}.  $\tens{T}(:,j,k)$, $\tens{T}(i,:,k)$ and $\tens{T}(i,j,:)$ denote vector slices of tensor $\tens{T}$. 
 
The operator  $\vec$ is the vectorization operator that stacks the columns of a matrix into a long column
vector, and $\Unvec$  is its reverse operator. $C^{k}$ is the set of functions having continuous $k$th derivatives. 
Finally, $\KK$ is either the real or the complex field.

\section{Description of Algorithms and Complexity Analysis }\label{sec: dcpt_alg}

This section presents the description of some CP decomposition algorithms. In order to support further results, the complexity per iteration is calculated here for each algorithm using  Landau's notation, denoted by $\mathcal{O}$, and counting only multiplications, as recommended in \cite{golub2012matrix}. 
From Section \ref{sec: als} to Section \ref{sec: cgrad}, we describe two classical algorithms known in literature: ALS and conjugate gradient \cite{comon2009tensor}. In Section \ref{sec: dcpd}, the DCPD algorithm is presented.

For the CP decomposition algorithms described in the following, the input parameter $R$ denotes the rank of the output tensor. Assuming  $R_{0}$  is the rank of the input tensor $\tens{T}$, if  $R_{0} \leq R$,   then the algorithms perform an exact decomposition. On the other hand, if  $R_{0} > R$, a lower rank-$R$ approximation is computed. 

\subsection{Alternating least squares (ALS)}\label{sec: als}
The most commonly used algorithm for solving the CP decomposition is  ALS  \cite{smilde2005multi}. The goal is to update alternately each factor matrix in each iteration by solving a least squares problem conditioned on previous updates of the other factor matrices. There is no guarantee of convergence to the global solution, nor even to a critical point. The implementation is quite simple and it is detailed  in Alg.\ref{ALS-alg}.

 \begin{algorithm}[ht!]
\SetKwInOut{Input}{input}\SetKwInOut{Output}{output}
 \Input{$\tens{T} \in \KK^{I_{1}\times I_{2}\times \cdots \times I_{N}}$: input data,\\
  R: rank parameter. }
 \Output{$\matr{A}^{(n)} \in  \KK^{I_{n}\times R}$, for  $n=1, \ldots, N$: factor matrices} 
 
 Initialize $\matr{A}^{(1)}, \ldots,  \matr{A}^{(N)}$
 
 \Repeat {\textit{some stopping  criterion is  satisfied}}{  \For{$n = 1$ to $N$}{
 $\matr{V}  \leftarrow  \matr{A}^{(1)\T}\matr{A}^{(1)} \hadam \cdots \hadam \matr{A}^{(n-1)\T}\matr{A}^{(n-1)} \hadam \matr{A}^{(n+1)\T}\matr{A}^{(n+1)} \hadam \cdots \hadam \matr{A}^{(N)\T}\matr{A}^{(N)} $
 
 $\matr{A}^{(n)} =  \matr{T}_{(n)} (\matr{A}^{(N)} \khatri \cdots \khatri \matr{A}^{(n+1)} \khatri \matr{A}^{(n-1)} \khatri \cdots \khatri \matr{A}^{(1)} )\matr{V}^{+}$
 
 }
 }
\caption{ALS algorithm}\label{ALS-alg}
\end{algorithm}

\bigskip

The complexity per iteration (\textit{repeat} loop) of ALS may be calculated as follows.
The computation of matrix $\matr{V}$ needs $(I_1+\ldots +I_{n-1} + I_{n+1} \ldots I_N)R^2 + (N-2)R^2$ operations (multiplications). The pseudo-inverse of $\matr{V}$ is calculated by resorting to an SVD. For a $m\times n$  rank-$r$ matrix, with $m \geq n \geq r$, the explicit calculation of diagonal, left singular, and right singular matrices require $2mn^2-2n^3/3, 5mr^2-r^3/3,$ and $5nr^2-r^3/3$ multiplications, respectively \cite{golub2012matrix}. Thus, assuming for simplicity that $\matr{V}$ is a non-singular matrix, the number of operations to calculate its pseudo inverse is $35R^3/3 + R^2$.  
For updating $\matr{A}^{(n)}$, we need 
$$
R  \!\!\prod_{i=1, i\neq n}^N \!\!\! I_{i}  +  R\prod _{j=1}^N I_j + I_nR^2
$$ 
multiplications. 
These calculations must be performed for each $n \in \{1,\ldots, N\}$. Thus, the number of operations per iteration for ALS is  dominated by the term  composed of the product of all dimensions. Hence,
\begin{equation}
\#op =  \mathcal{O} \bigg\{NR\prod \limits_{j=1}^N I_j \bigg\}.
\end{equation}

\subsection{Conjugate gradient (CG) }\label{sec: cgrad}

The conjugate gradient algorithm (CG)  is  a faster algorithm than the well known gradient descent \cite{comon2009tensor}. Here, we use the optimum step size and the Polak-Ribi\`{e}re implementation \cite{polak1997optimization} for updating the parameter $\beta$ in the algorithm presented in Alg.\ref{CG-alg}. 
The number of operations for computing each vector $\vect{g}_{n}$ is 
$$
R^2 \sum \limits_{j=1}^{N}I_j  + (N-2)R^2 +   I_nR^2 + R  \!\!\prod\limits_{i=1, i\neq n}^N \!\!\! I_{i}  +  R\prod \limits_{j=1}^{N}I_j
$$ 
The computation of the step size $\mu^{\star}$ is dominated by the number of multiplications needed to determine all the coefficients but one of a $2N$-degree polynomial generated from the enhanced line search (ELS) method, which is given by \cite{rajih2008enhanced}: 
$$
\left(2^{N}R+ \mathcal{O} \{N^{2}\} \right)\prod \limits_{j=1}^{N}I_j  
$$

 \begin{algorithm}[ht!]
\SetKwInOut{Input}{input}\SetKwInOut{Output}{output}
 \Input{$\tens{T} \in \KK^{I_{1}\times I_{2}\times \cdots \times I_{N}}$: input data,\\
  R: rank parameter. }
 \Output{$\matr{A}^{(n)} \in  \KK^{I_{n}\times R}$, for  $n=1, \ldots, N$: factor matrices} 
 
 Initialize $\matr{A}^{(1)}, \ldots,  \matr{A}^{(N)}$; 
 
$\vect{p} \leftarrow [vec^{\T}(\matr{A}^{(1)}) \cdots vec^{\T}(\matr{A}^{(N)})]^{\T}$;

 \Repeat {\textit{some stopping  criterion is  satisfied}}{ 
 
  \For{$n = 1$ to $N$}{
  
 Compute the gradient $\vect{g}_{(n)}$ with respect $vec(\matr{A}^{(n)})$;
  
  }
  
  $\vect{g} \leftarrow [\vect{g}^{\T}_{(1)} \ldots \vect{g}^{\T}_{(N)}]^{\T}$;
  
 \If {$n=1$}{ $\vect{d} \leftarrow -\vect{g}$;} 
 
  Compute the optimum step size $\mu^{\star}$;
  
 Update $\beta$ according to Polak-Ribi\`{e}re;
  
  $\Delta \vect{p} \leftarrow \mu^{\star} \vect{d}$;
   
   $\vect{p}  \leftarrow \vect{p}  + \Delta \vect{p}$;
   
     $\vect{d} \leftarrow - \vect{g} + \beta \vect{d}$;
 
  }
  
  Extract $\matr{A}^{(n)}$ from $\vect{p}$, for $n = 1, \dots, N.$
  
\caption{CG algorithm}\label{CG-alg}
\end{algorithm}

The parameter $\beta$ requires only $1+ 2R\sum_{j=1}^{N}I_j$ multiplications. Hence, the CG algorithm with ELS has a total complexity given by

 \begin{equation}
\#op =  \mathcal{O} \bigg\{\left( (2^{N}+1)R  +  N^{2}\right)\prod \limits_{j=1}^{N}I_j \bigg\}.
\end{equation}

\subsection{Deflation-based CP decomposition (DCPD) }\label{sec: dcpd}

The computation of a rank-1 approximation is the key of the DCPD algorithm. We present here two methods for computing a rank-1 approximation referred to as THOSVD and SeROAP \cite{AlexIcassp2015}.

\subsubsection{ Truncation of higher order singular value decomposition - THOSVD}

The algorithm is described in Alg.\ref{HOSVD-alg}.
For computing the first right singular vector, we do not need compute the complete SVD. According to \cite{comon1990tracking}, we can compute the best rank-1 approximation  of a $m\times n$ matrix in $k$ steps by using the Lanczos algorithm, with a complexity  $\mathcal{O} \{2kmn \}$.  Hence, the accumulated complexity computed for all $\vect{u}_{n}$ is equal to $\mathcal{O} \{2 Nk  \prod_{j=1}^{N}I_j \}$. 
The computation of $\tens{U}$ requires $\prod_{j=1}^{N}I_j$ flops. The contraction to obtain $\lambda$ also needs $\prod_{j=1}^{N}I_j$ operations.
To sum up, the total number of operations of THOSVD is of order:
$$
\mathcal{O}\bigg\{(2 Nk+2)\prod_{j=1}^N I_j\bigg\} 
$$ 

 \begin{algorithm}[ht!]
\SetKwInOut{Input}{input}\SetKwInOut{Output}{output}
 \Input{$\tens{T} \in \KK^{I_{1}\times I_{2}\times \cdots \times I_{N}}$: input data}
 \Output{ $\tens{X}  \in  \KK^{I_{1}\times I_{2}\times \cdots \times I_{N}}$: rank-1 approximation} 
 
 \For{$n = 1$ to $N$}{

$\vect{u}_{n} \leftarrow$ first left singular vector of  $\matr{T}_{(n)}$ ;
  
 }
$\tens{U} \leftarrow \out_{n=1}^{N}\vect{u}_{n}$;

$\lambda \leftarrow \langle \tens{T}, \tens{U} \rangle$;

$\tens{X}  \leftarrow \lambda \cdot \tens{U}$.

\caption{THOSVD algorithm}\label{HOSVD-alg}
\end{algorithm}

\subsubsection{Sequential rank-1 approximation and projection - SeROAP}

 Without loss of generality, consider $I_1 \geq I_2 \geq \ldots \geq I_N$. The SeROAP algorithm \cite{AlexIcassp2015} goes along the lines depicted in Alg.\ref{SeROAP-alg}. In the first \textit{for} loop, we compute $N-2$ right singular vectors of matrices  whose size is successively reduced. For this step, the complexity is $\mathcal{O}\{ 2k\sum_{i=1}^{N-2} (\prod_{j=i}^{N}I_j) \}$.  The computation of vectors $\vect{u},\vect{v}$ and $\vect{w}$ has complexity  $\mathcal{O}\{(2k+1)I_{N-1}I_{N} \}$. 

Next, the second \textit{for} loop  performs $N-2$ successive projections of the rows of matrices $\matr{V}_{n}$ onto the vectors $\vect{w}$. We need here $2\sum_{i=1}^{N-2}(\prod_{j=i}^{N}I_j) $ operations.

\begin{algorithm}[ht!]
\SetKwInOut{Input}{input}\SetKwInOut{Output}{output}
 \Input{$\tens{T} \in \KK^{I_{1}\times I_{2}\times \cdots \times I_{N}}$: input data}
 \Output{$\tens{X} \in  \KK^{I_{1}\times I_{2}\times \cdots \times I_{N}}$: rank-1 approximation}

 $\matr{V}_{0} \leftarrow  \matr{T}_{(1)}$;
 
 $\matr{V} \leftarrow \matr{V}_{0}$;
   
 \For{$n = 1$ to $N-2$}{

$\vect{v}_{n} \leftarrow$ first right singular vector of $\matr{V}$;

$\matr{V}_{n} \leftarrow  Unvec (\vect{v}_{n}) \in \KK^{I_{n+1}\times I_{n+2}  I_{n+3}\cdots I_{N}}$;

$\matr{V} \leftarrow \matr{V}_{n}$;
  
 }
 $ (\vect{u},\vect{v}) \leftarrow$ first left and right singular vectors of $\matr{V}$;
 $ \vect{w} \leftarrow  \vect{v}^{*}\kron \vect{u}$; 
  
 \For{$n = N-2$ to $1$}{
 $\matr{X}_{(n)} \leftarrow  \left(\matr{V}_{n-1} \vect{w}\right)\vect{w}^{H}$; 
 
 $\vect{w} =  vec({\matr{X}_{(n)}}) $;
 }
 $\matr{X}_{(1)}$ is the mode-$1$ unfolding of  $\tens{X}$.

\caption{SeROAP algorithm}\label{SeROAP-alg}
\end{algorithm}

For large dimensions and small $N$, the complexity of SeROAP is dominated by: 
$$
\mathcal{O} \bigg\{(2k+2)\prod \limits_{j=1}^{N}I_j \bigg\}
$$ 
which can be significantly smaller than that of THOSVD.
An example of typical execution\footnote{Matlab codes are available at  http://www.gipsa-lab.grenoble-inp.fr/$\sim$pierre. comon/TensorPackage/tensorPackage.html.} is given in the Appendix.  

\subsubsection{ Description and complexity of DCPD}
 
The DCPD is an iterative deflation algorithm \cite{AlexIcassp2015} that computes the CP decomposition for real or complex tensors.  As summarized in Alg.\ref{DCPD-alg}, it proceeds as follows.  In the first \textit{for} loop, we compute the rank-1 tensors $\tens{X}[1,1], \ldots, \tens{X}[R,1]$  by successive rank-1 approximations and subtractions. Since the rank of the tensor does not decrease with subtractions in general \cite{StegC10:laa}, a residual  $\tens{E}[R,1]$  is then produced. 
In the iterative process (\textit{repeat} loop), a new rank-1 component is generated from the sum of the previous residual and  $\tens{X}[1,1]$,  and a new residual  $\tens{E}[1,2]$  is produced with the subtraction  $\tens{Y}[1,2] - \tens{X}[1,1]$.  The tensor $\tens{Y}[1,2]$ is updated within the \textit{if-else} condition.  By applying the same procedure to the other components, we update all $R$ rank-1 tensors, so that another residual  $\tens{E}[R,2]$  is generated in the end of the second \textit{for} loop. The second loop continues to execute until some stopping criterion is satisfied,  and all rank-1 components of $\tens{T}$ can be recovered.

\begin{algorithm}[ht!]
\SetKwInOut{Input}{input}\SetKwInOut{Output}{output}
\Input{$\tens{T} \in \KK^{I_{1}\times I_{2}\times \cdots \times I_{N}}$: input data,\\
  R: rank parameter. \\
  $\phi$: an algorithm computing a rank-1 \\approximation}
 \Output{$\tens{X}_{r} \in \KK^{I_{1}\times I_{2}\times \cdots \times I_{N}}$, for  $r=1, \ldots, R$: rank-1 components} 
   
  $\tens{Y}[1,1] \leftarrow \tens{T}$;
  
  \For{$r = 1$ to $R$}{
  $ \tens{X}[r,1] \leftarrow \phi (\tens{Y}[r,1])$;
    
   \eIf {$r<R$}{ $\tens{Y}[r+1,1] \leftarrow \tens{Y}[r,1]-\tens{X}[r,1];$} 
      {$\tens{E}[R,1] \leftarrow  \tens{Y}[R,1]-\tens{X}[R,1];$ }}

$l \leftarrow 2;$\\
 \Repeat {\textit{some stopping criterion is  satisfied}}{  \For{$r = 1$ to $R$}{
    
    \eIf {$r >1$}{
    $\tens{Y}[r,l] \leftarrow \tens{X}[r,l-1] + \tens{E}[r-1,l] ;$ 
   }
   { $\tens{Y}[1,l] \leftarrow \tens{X}[1,l-1] + \tens{E}[R,l-1];$ }
   $\tens{X}[r,l] \leftarrow \phi (\tens{Y}[r,l])$;\\
   $\tens{E}[r,l] \leftarrow \tens{Y}[r,l] - \tens{X}[r,l] $;\\
  }
  $l \leftarrow l + 1;$\\
 }
 \
\ForEach {$r \in [1,\ldots,R]$}{ $\tens{X}_r \leftarrow \tens{X}[r,l];$}
  \caption{DCPD algorithm}\label{DCPD-alg}
\end{algorithm}

The complexity per iteration is dominated by the rank-1 approximation function $\phi$, which is computed $R$ times. Therefore, the complexity of DCPD is 
 \begin{equation}
\#op =  \mathcal{O} \left\{(2Nk+2)R\prod \limits_{j=1}^{N}I_j \right\}, 
\end{equation} for the T-HSOVD algorithm, and

\begin{equation}
\#op =  \mathcal{O} \left\{(2k+2)R\prod \limits_{j=1}^{N}I_j \right\}, 
\end{equation} for the SeROAP algorithm.

\section{rank-1 Approximation}\label{sec:rank_one_res} 
This section is divided in two parts.  The first one presents a more detailed study of the THOSVD and SeROP rank-1 approximation methods. For three-way tensors, we  show  that SeROAP is a better choice than THOSVD because the  former presents a better rank-1 approximation, which ensures a more probable monotonic  decrease of the residual  $\|\tens{E}[R,l]\|$ within the  DCPD algorithm \cite{AlexIcassp2015}. A new rank-1 approximation algorithm  is  described in a second part, and is proved to perform better than any other rank-1 approximation  method.

\subsection{THOSVD vs SeROAP} \label{sec: THvsSeR}

Since we do not have at our disposal an efficient method to compute quickly the best rank-1 approximation of a tensor, we should compute a suboptimal rank-1 approximation in a tractable way. The THOSVD and SeROAP algorithms can perform this task, as presented in Section \ref{sec: dcpt_alg}. The question that arises is: which algorithm performs the best? The proposition  \ref{prop:THSe} below  shows  that SeROAP  performs better than THOSVD for three-way tensors.  For simplicity, the notation of the unfolding matrices does not present indices in  this section.

\smallskip
\begin{proposition}\label{prop:THSe}
Let $\tens{T} \in \KK^{I_{1}\times I_{2}\times I_{3}}$ be a $3$-order tensor. Let also $\phi^{\textrm{TH}}(\tens{T})$ and  $\phi^{\textrm{Se}}(\tens{T})$ be the  rank-1 approximations delivered by THOSVD and SeROAP algorithms, respectively. 
Then the inequality  $\| \tens{T} -\phi^{\textrm{Se}}(\tens{T})\| \leq \| \tens{T} -\phi^{\textrm{TH}}(\tens{T})\|$ holds.
\end{proposition}

\begin{proof}
Let $\matr{T}$, $\matr{T}_{\phi}^{\textrm{TH}}$  and $\matr{T}_{\phi}^{\textrm{Se}}$ be some mode unfolding of tensors $\tens{T}$, $\phi^{\textrm{TH}}(\tens{T})$ and  $\phi^{\textrm{Se}}(\tens{T})$, respectively. Assuming mode-1 unfolding for the THOSVD algorithm, we have 
$$
 \| \matr{T} -\matr{T}_{\phi}^{\textrm{TH}}\|^{2}=  \| \matr{T} - \lambda \vect{u}_{1}(\vect{u}_{3} \kron \vect{u}_{2})^{T}\|^{2},
$$
where $\lambda, \vect{u}_{1}, \vect{u}_{2}$, and $\vect{u}_{3}$ are obtained from Alg. \ref{HOSVD-alg}. Since $\lambda$ is the contraction of $\tens{T}$ on $\tens{U}$, we plug it into the previous equation and we obtain, after simplifications, 
$$
\| \matr{T} -\matr{T}_{\phi}^{\textrm{TH}}\|^{2}= \| \matr{T} \|^{2} -|\lambda|^2,
$$ 
with $\lambda = \vect{u}_1^H\matr{T}(\vect{u}_3^*\kron\vect{u}_2^*)$.
Yet, $\vect{u}_1^H\matr{T}=\|\matr{T}\|_2\,\vect{v}_1^H$ since $(\vect{u}_1, \vect{v}_1,\, \|\matr{T}\|_2)$ is the dominant singular triplet of matrix $\matr{T}$. 
Hence $|\lambda|^2=\|\matr{T}\|_2^2 \, |\vect{v}_{1}^{H}(\vect{u}_{3}^{*} \kron \vect{u}_{2}^{*} )|^2$.

On the other hand for SeROAP, we have 
$$ 
\| \matr{T} -\matr{T}_{\phi}^{\textrm{Se}}\|^{2}= \| \matr{T} -\matr{T}\vect{w}\vect{w}^{H}\|^{2}= 
$$ 
$$
\|\matr{T}\|_F^{2} + \|\matr{T}\vect{w}\|_2^2 \, \|\vect{w}\|_2^2-2\|\matr{T}\vect{w}\|_2^2
= \| \matr{T} \|^{2} - \vect{w}^{H} \matr{T}^{H}\matr{T}\vect{w},
$$ 
 where $\vect{w}$ is the same vector computed before the second \textit{for} loop given in Alg. \ref{SeROAP-alg} for $3$-order tensors. 
The eigenvalue decomposition of $ \matr{T}^{H}\matr{T}$ can be expressed by $$  \matr{T}^{H}\matr{T} = \|\matr{T}\|_{2}^{2}\vect{v}_{1}\vect{v}_{1}^{H} + \matr{S},$$ where $ \matr{S}$ is a semidefinite positive matrix. 
Hence, we have
$$
\vect{w}^{H} \matr{T}^{H}\matr{T}\vect{w} = \|\matr{T}\|_{2}^{2}\vect{w}^{H} \vect{v}_{1}\vect{v}_{1}^{H}\vect{w} + c, 
$$ 
with $c \geq 0$.
To complete the proof of the proposition, we just need to show that $|\vect{v}_{1}^{H}\vect{w}|^2 \geq |\vect{v}_{1}^{H}(\vect{u}_{3}^{*} \kron \vect{u}_{2}^{*} )|^2$, 
or equivalently that
$$
 |\langle\vect{w},\vect{v}_{1} \rangle| \geq |\langle \vect{u}_{3}^{*} \kron \vect{u}_{2}^{*} ,\vect{v}_{1}\rangle|.
$$ 
This is true, because $\vect{w}$ is by construction (cf. Alg. \ref{SeROAP-alg}) the vector closest to $\vect{v}_{1}$ among all vectors of the form $\vect{a}\kron\vect{b}$ where $\vect{a}$ and $\vect{b}$ have unit norm. 
\end{proof}
\smallskip

\subsection{Coupled-eigenvalue rank-1 approximation} \label{sec:CERankOne}

This section presents an alternating eigenvalue method for three-way tensors that can improve local solutions obtained from any other rank-1 approximation method (e.g. SeROAP and THOSVD algorithms). Actually, simulations have shown that the global solution is always attained if the initial approximation is  close enough.

Let $\vect{t}_{i_3}$ be the vectorization of slice $i_3$, $ 1\leq i_3 \leq I_3$ (we have chosen the third mode) of tensor $\tens{T}$. The rank-1 approximation problem can be stated as
\begin{equation}\label{rank1:eq}
\begin{array}{c}
[\vect{\alpha}_{opt},\vect{x}_{opt},\vect{y}_{opt}] = arg\displaystyle\min_{\vect{\alpha},\vect{x},\vect{y}}  \Upsilon(\vect{x},\vect{y},\vect{\alpha})  \\
\\\vspace{-0.7cm}\\
s.t. ~ \|\vect{x}\|_2 =1, \, \|\vect{y}\|_2 =1,
\end{array} 
\end{equation} with 
 $\Upsilon(\vect{x},\vect{y},\vect{\alpha})=\sum_{i_{3}=1}^{I_3}\|\vect{t}_{i_3}-\alpha_{i_3}(\vect{y}^{*}\kron \vect{x})\|_{2}^{2}$ and $\vect{\alpha} = [\alpha_1 \cdots \alpha_{I_3}]$. 

Plugging the optimal value of $\alpha_{i_3}$ into Problem (\ref{rank1:eq}),  namely  $(\vect{y}^{*}\kron \vect{x})^H\vect{t}_{i_3}$,  we can rewrite it as the following equivalent maximization problem
\begin{equation}\label{rank1:eq_max}
\begin{array}{c}
[\vect{x}_{opt},\vect{y}_{opt}] = arg\displaystyle\max_{\vect{x},\vect{y}} \vect{z}^{H}\matr{M}\vect{z} \\
\\\vspace{-0.7cm}\\
s.t. \hspace{0.1cm} ~ \vect{z} = \vect{y}^{*}\kron \vect{x},
  \, \|\vect{y}\|_2 =1, \|\vect{x}\|_2 =1,
\end{array} 
\end{equation} where $\matr{M} = \sum_{i_{3}=1}^{I_3}\vect{t}_{i_3}^{} \vect{t}_{i_3}^{H}$.

Now, we decompose $\matr{M}$ as a sum of Kronecker products. This can be done by reshaping $\matr{M}$ and applying the SVD decomposition \cite{van1993approximation}.
Thus, $\matr{M}$ can be given by 
\[ \matr{M} = \sum \limits_{r=1}^{R'}  \matr{Q}^{(r)}   \kron \matr{P}^{(r)},
\]
with the  Hermitian  matrices $\matr{P}^{(r)} \in \mathbb{K}^{I_{1}\times I_{1}}$ and $\matr{Q}^{(r)} \in \mathbb{K}^{I_{2}\times I_{2}}$. $R'$ is the \emph{Kronecker rank} of $\matr{M}$ satisfying $R' \leq \min(I_{1}^{2},I_{2}^{2})$. Substituting $\matr{M}$ into Problem (\ref{rank1:eq_max}), we have:
\begin{equation}\label{rank1:eigen}
\begin{array}{c}
[\vect{x}_{opt},\vect{y}_{opt}] = arg\displaystyle\max_{\vect{x},\vect{y}} \sum \limits_{r=1}^{R'} (\vect{y}^{H}\matr{Q}^{(r)*}\vect{y} )(\vect{x}^{H}\matr{P}^{(r)}\vect{x} )\\
\\\vspace{-0.7cm}\\
s.t.  \|\vect{y}\|_2 =1,   \, \|\vect{x}\|_2 =1.
\end{array} 
\end{equation}  
Let $\mathcal{L}$ be the Lagrangian function given by
\small
\[ \mathcal{L} = -\sum \limits_{r=1}^{R'} (\vect{y}^{H}\matr{Q}^{(r)*}\vect{y} )(\vect{x}^{H}\matr{P}^{(r)}\vect{x} ) + \eta_{1} (\|\vect{y}\|_{2}^{2} -1)+ \eta_{2} (\|\vect{x}\|_{2}^{2} -1), \]
\normalsize
where  $\eta_{1}$ and $\eta_{2}$ are the Lagrange multipliers.  By computing the critical points, we obtain a pair of coupled eigenvalue problems
 
\begin{equation}\label{CEy}
\left [ \begin{array}{ccc}
\vect{y}^{H}\matr{A}_{(1,1)}\vect{y} & \cdots & \vect{y}^{H}\matr{A}_{(1,I_{1})}\vect{y} \\
\vdots & \ddots & \vdots \\
\vect{y}^{H}\matr{A}_{(I_{1},1)}\vect{y} & \cdots &  \vect{y}^{H}\matr{A}_{(I_{1},I_{1})}\vect{y}  \\
\end{array} \right]\vect{x} = \lambda \vect{x}
\end{equation} and
\begin{equation}\label{CEx}
\left [ \begin{array}{ccc}
\vect{x}^{H}\matr{B}_{(1,1)}\vect{x} & \cdots & \vect{x}^{H}\matr{B}_{(1,I_{2})}\vect{x} \\
\vdots & \ddots & \vdots \\
\vect{x}^{H}\matr{B}_{(I_{2},1)}\vect{x} & \cdots &  \vect{x}^{H}\matr{B}_{(I_{2},I_{2})}\vect{x}  \\
\end{array} \right]\vect{y} = \lambda \vect{y},
\end{equation} where  $\lambda = \eta_1 = \eta_2$, $\matr{A}_{(m,n)} = \sum \limits_{r=1}^{R'}P_{mn}^{(r)}\matr{Q}^{(r)*}$, and  $\matr{B}_{(k,l)} = \sum \limits_{r=1}^{R'}Q_{kl}^{(r)*}\matr{P}^{(r)}$, with  $1 \leq m,n \leq I_1$ and $1 \leq k,l \leq I_2$.  

\smallskip

The coupled-eigenvalue  algorithm is presented in Alg. \ref{CoupledEig-alg}.   We can initialize the algorithm by computing $\vect{x}_{0}$ and $\vect{y}_{0}$ from the rank-1 approximation solution obtained with SeROAP, THOSVD or any other rank-1 approximation method.

\begin{algorithm}[ht!]
\SetKwInOut{Input}{input}\SetKwInOut{Output}{output}
 \Input{$\phi(\tens{T})$: rank-1 approximation}
 \Output{$\phi^{\star}(\tens{T})$: improved rank-1 approximation} 

Compute $\vect{x}_{0}$  from $\phi(\tens{T})$ as \\
\small  $\vect{x}_{0} \leftarrow \phi(\tens{T})(:,i_2,i_3)/\| \phi(\tens{T})(:,i_2,i_3)\|_2 \hspace{1mm} \text{for some} \hspace{1mm}  i_{2}, i_{3}$

 $t \leftarrow 0$\\
  \Repeat {  \textit{some stopping  criterion is  satisfied}  }{ 
  
  Set $\vect{x}=\vect{x}_{t}$  in eigenvalue problem (9)  and take $\vect{y}_{t+1}$  as the  eigenvector whose eigenvalue is maximum;\\
  Set $\vect{y}=\vect{y}_{t+1}$  in eigenvalue problem (8)  and take  $\vect{x}_{t+1}$  as the  eigenvector whose eigenvalue is maximum;  \\
  $t \leftarrow t + 1$
  } 
 \For{$i_3 = 1$ to $I_3$}{
  $\alpha_{i_3}^{\star} \leftarrow \left \langle \vect{t}_{i_3}, \vect{y}_{t}^{*} \kron \vect{x}_{t} \right \rangle$;
   
 }
 $\phi^{\star}(\tens{T}) \leftarrow \vect{x}_{t} \otimes \vect{y}_{t} \otimes \vect{\alpha}^{\star}$;
 
\caption{CE rank-1 approximation. Above, we chose to start with $\vect{x}$, but we could equivalently have started with $\vect{y}$. }\label{CoupledEig-alg}
\end{algorithm}

The complexity per iteration of the CE algorithm is dominated by the construction of the matrices in the  LHS of (\ref{CEy}) and (\ref{CEx}), which is  of order  $\mathcal{O}\{\min(I_{1}^{2}I_{2}^2, I_{1}^{2}I_{3}^2,I_{2}^{2}I_{3}^2)\}$. Suppose $I_3$ is the largest dimension. If $I_{3} \gg I_{1}I_{2}$, then we can take advantage of the CE algorithm in terms of complexity in comparison with the ALS algorithm. Indeed, the complexity per iteration of ALS for rank-1 approximation is of order $\mathcal{O}(3I_{1}I_{2}I_{3})$, which is higher than that of the CE algorithm in this case.  Notice, however,  that a properly comparison makes sense if the same initialization is employed in both algorithms.

The following proposition shows that the above algorithm improves (in worst case the solution remains the same) any rank-1 approximation algorithm.
\smallskip
\begin{proposition}\label{prop:eigen}
Let $\phi(\tens{T})$ be a rank-1 approximation of a three-way tensor $\tens{T}$. If $\phi(\tens{T})$ is the input of the CE algorithm and $\phi^{\star}(\tens{T})$ the output, then the  inequality  $\| \tens{T} -\phi^{\star}(\tens{T})\| \leq \| \tens{T} -\phi(\tens{T})\|$ holds.
\end{proposition}
\smallskip

\begin{proof} 
Plugging the expression of $\matr{A}_{(m,n)}$ into equation (\ref{CEy}), we obtain, after simplifications,
$$
 \lambda = \sum \limits_{r=1}^{R'} (\vect{y}^{H}\matr{Q}^{(r)*}\vect{y} )(\vect{x}^{H}\matr{P}^{(r)}\vect{x} ), 
$$ 
which is the objective function of  Problem (\ref{rank1:eigen}). The same result is obtained when the matrix $\matr{B}_{(k,l)}$ is plugged into equation (\ref{CEx}).
Now 
$\forall t \geq 1$, let $\lambda_{t}^{(\vect{x})}$ and  $\lambda_{t}^{(\vect{y})}$  be the  maximal eigenvalues  whose eigenvectors are $\vect{x}_t$  and $\vect{y}_t$, respectively. 

The eigenpair $(\lambda_{t+1}^{(\vect{y})}, \vect{y}_{t+1})$  obtained by solving  equation (\ref{CEx}) with $\vect{x}=\vect{x}_{t}$, is solution of the maximization problem
$$
\lambda_{t+1}^{(\vect{y})}= \displaystyle\max_{\|\vect{y}\|_{2}=1}  \sum \limits_{r=1}^{R'} (\vect{y}^{H}\matr{Q}^{(r)*}\vect{y} )(\vect{x}_{t}^{H}\matr{P}^{(r)}\vect{x}_{t} ).
$$
Also, the eigenpair $(\lambda_{t+1}^{(\vect{x})}, \vect{x}_{t+1})$ obtained by setting $\vect{y}=\vect{y}_{t+1}$ in  equation (\ref{CEy}), is solution of the problem
$$
\lambda_{t+1}^{(\vect{x})} = \displaystyle\max_{\|\vect{x}\|_{2}=1}  \sum \limits_{r=1}^{R'} (\vect{y}_{t+1}^{H}\matr{Q}^{(r)*}\vect{y}_{t+1} )(\vect{x}^{H}\matr{P}^{(r)}\vect{x} ).
$$

Since 
\begin{align*}
 \displaystyle\max_{\|\vect{x}\|_{2}=1}  \sum \limits_{r=1}^{R'} (\vect{y}_{t+1}^{H}&\matr{Q}^{(r)*}\vect{y}_{t+1} )(\vect{x}^{H}\matr{P}^{(r)}\vect{x} ) &\\ &= \sum \limits_{r=1}^{R'} (\vect{y}_{t+1}^{H}\matr{Q}^{(r)*}\vect{y}_{t+1} )(\vect{x}_{t+1}^{H}\matr{P}^{(r)}\vect{x}_{t+1} ),
\end{align*} it follows in particular that 

\begin{align*}
 \sum \limits_{r=1}^{R'} (\vect{y}_{t+1}^{H}&\matr{Q}^{(r)*}\vect{y}_{t+1} )(\vect{x}_{t}^{H}\matr{P}^{(r)}\vect{x}_{t} )  &\\ &\leq \sum \limits_{r=1}^{R'} (\vect{y}_{t+1}^{H}\matr{Q}^{(r)*}\vect{y}_{t+1} )(\vect{x}_{t+1}^{H}\matr{P}^{(r)}\vect{x}_{t+1} ),
\end{align*} 
which implies that $\lambda_{t+1}^{(\vect{y})} \leq \lambda_{t+1}^{(\vect{x})}$.

Similarly, plugging $\vect{x}_{t+1}$ into equation (\ref{CEx}), we can conclude that $\lambda_{t+1}^{(\vect{x})} \leq  \lambda_{t+2}^{(\vect{y})}$ for the reason that 

\begin{align*}
 \sum \limits_{r=1}^{R'} (\vect{y}_{t+1}^{H}&\matr{Q}^{(r)*}\vect{y}_{t+1} )(\vect{x}_{t+1}^{H}\matr{P}^{(r)}\vect{x}_{t+1} ) &\\ & \leq \sum \limits_{r=1}^{R'} (\vect{y}_{t+2}^{H}\matr{Q}^{(r)*}\vect{y}_{t+2} )(\vect{x}_{t+1}^{H}\matr{P}^{(r)}\vect{x}_{t+1} ).
\end{align*}

Hence, the sequence 
$$
 \{\Upsilon_t\}_{t \in \mathbb{N}}  = \{\ldots, \lambda_{t}^{(\vect{y})},\lambda_{t+1}^{(\vect{x})},\lambda_{t+1}^{(\vect{y})},\lambda_{t+2}^{(\vect{x})},\ldots \}
$$ 
is monotonically non-decreasing. The same conclusion would be achieved if we begin by plugging  $\vect{x}_t$ into equation (\ref{CEx}).

Now, let $\phi(\tens{T}) = \vect{x}_{0} \otimes \vect{y}_{0} \otimes \vect{\alpha}_{0}$ be a rank-1 approximation obtained with any  other method. Assume $\vect{x}_{0}$ and $\vect{y}_{0}$  are  unit vectors, and define $\lambda_{0} = \sum \limits_{r=1}^{R'} (\vect{y}_{0}^{H}\matr{Q}^{(r)*}\vect{y}_{0} )(\vect{x}_{0}^{H}\matr{P}^{(r)}\vect{x}_{0} )$. By setting $\vect{x}=\vect{x}_{0}$ in equation (\ref{CEx}) in the first iteration  (a similar operation would be possible for $\vect{y}_{0}$ in equation (\ref{CEy})),  we clearly have $\lambda_{0} \leq  \lambda_{1}^{(\vect{y})} \leq \lambda_{t_{max}}^{(\vect{y})}$, where $t_{max}$ is the iteration in which the stopping criterion is satisfied.

Since the optimization problems (\ref{rank1:eq}) and (\ref{rank1:eigen}) are equivalent,   $\alpha_{i_3}^{\star}, 1\leq i_3 \leq I_3$, can be obtained by performing the scalar product between  vectors $ \vect{t}_{i_3}$ and  $\vect{y}_{t_{max}}^{*} \!\kron \vect{x}_{t_{max}}$ (which is equivalent to contracting tensor $\tens{T}$ on $\vect{x}_{t_{max}}$ and $\vect{y}_{t_{max}}^{*}$).  Hence, the tensor $$\phi^{\star}(\tens{T}) = \vect{x}_{t_{max}} \otimes \vect{y}_{t_{max}} \otimes \vect{\alpha}^{\star}$$ is a better rank-1 approximation  of $\tens{T}$ than $\phi(\tens{T})$,  implying  $\| \tens{T} -\phi^{\star}(\tens{T})\| \leq \| \tens{T} -\phi(\tens{T})\|$.
\end{proof} 

\medskip

\begin{proposition}\label{prop:CEconv} 
 For any input $\phi(\tens{T})$, the CE algorithm converges to a stationary point.
\end{proposition}
\smallskip
\begin{proof} 
In the proof of Proposition \ref{prop:eigen}, we have shown that $\{\Upsilon_t\}$  is monotonically non-decreasing for any input $\phi(\tens{T})$. Let $p^{\star}$ be the maximum of the objective (\ref{rank1:eigen}). Since the best rank-1 approximation problem  always has a solution, then  $p^{\star} < \infty$.  But $\max_{x}\Upsilon(\vect{x},\vect{y_{t+1}},\vect{\alpha}) \le \max_{x,y}\Upsilon(\vect{x},\vect{y},\vect{\alpha})$, which 
 implies that $\{\Upsilon_t\}_{t \in \mathbb{N}}$ is bounded above by $ p^{\star}$. 
 Since $\{\Upsilon_t\}$ is a real non decreasing sequence bounded above, it converges to a limit $\Upsilon^\star$, $\Upsilon^\star\le p^\star$.
\end{proof}
\smallskip

\section{Deflation}\label{sec:pfDCPD}

In \cite{AlexIcassp2015}, we proved that for a rank-$R$ tensor, the normalized residual $(\|\tens{E}[R,l]\|)_{l \in \mathbb{N}_{> 0}}$ is a monotonically decreasing sequence when the best rank-1 approximation is assumed within DCPD. In this section, a thorough theoretical analysis and new results are presented. Based on a geometric approach, we sketch an analysis of the convergence of the DCDP algorithm, including a conjecture that it converges to an exact decomposition with high probability when tensors within $\set{T}^{(R)} = \{\tens{T} \in \set{T}: ~ \rank{ \tens{T}} \leq R\}$ are distributed according to  absolutely continuous probability measures.

\begin{figure}
 \centering{


\begin{tikzpicture}

%

\draw[line width= 1pt, draw=black] (0,0) circle (2);

\draw[line width= 1pt, draw=blue,-triangle 90,fill=blue]  (-4.5,0) -- (0,0);


\draw [color = gray] (0,0) --(70:2);
\draw [color = gray] (0,0) --(110:2);
\draw [color = gray] (0,0) --(-70:2);
\draw [color = gray] (0,0) --(-110:2);
\path [fill=gray] (0,0) -- (110:2) arc (110:70:2);
\path [fill=gray] (0,0) -- (-110:2) arc (-110:-70:2);

\draw [line width= 1pt, style=dashed, color = black] (0,0) --(3,0) coordinate (x axis);


\draw[line width= 1pt, draw=red,-triangle 90,fill=red]  (0,0) -- (55:3.5);

\draw[line width= 1pt, draw=gray,-triangle 90,fill=gray]  (-4.5,0) -- (55:3.5);

\draw[very thick,black]
(-30:0.9cm) -- node[left=1pt] {$\beta$}  (-30:0.9cm);

\draw[very thick,black]
(15:1.6cm) -- node[left=1pt] {$\gamma[r,l]$}  (15:1.6cm);

\draw [line width= 1pt, style=dashed, color = black] (0,-3) --(0,3) coordinate (y axis);

\draw [draw=black](0:0.6) arc (0:55:0.6);

\draw [draw=black](0:0.4) arc (0:-70:0.4);
\draw [draw=black] (0:0.5) arc (0:-70:0.5);
\draw (-3,-0.26) node { $\tens{X}[r,l-1]$};
\draw [draw=gray, text=gray] (-4.5,0) -- (55:3.5) node[midway,sloped,above]{ $\tens{Y}[r,l]$};
\draw[text=red] (40: 2.8) node {$\tens{E}[r-1,l]$};

\end{tikzpicture}
}
   \vspace{-2.0ex}
 \caption{Visualization of the residual in an $n$-sphere for some iteration $l$ of DCPD algorithm. }\label{fig:cone}
\end{figure}
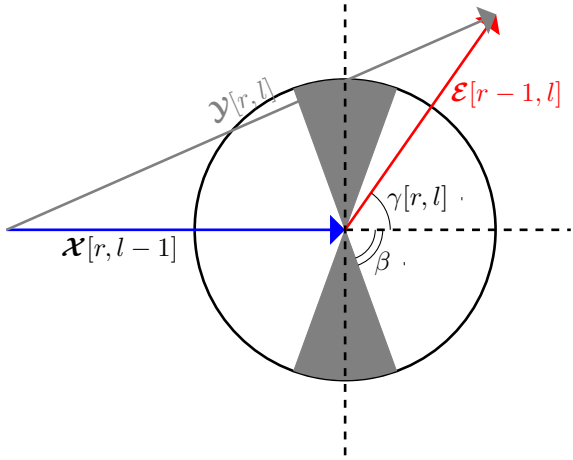

First, let us take a closer look at the 2D geometric interpretation of the DCPD algorithm. Figure \ref{fig:cone}  depicts  the $[r,l]$-iteration  for $r > 1$, so that $\gamma[r,l]$ is the  angle between the  tensors  $\tens{E}[r-1,l]$ and  $\tens{X}[1,l-1]$.  
For $r=1$, the residual  $\tens{E}[r-1,l]$  can be just replaced with $\tens{E}[R,l-1]$ in the figure,  and $\gamma[1,l]$ is  then  defined from  $\tens{E}[R,l-1]$ and $\tens{X}[1,l-1]$. 

Before stating some theoretical results on the DCPD algorithm, we present a fundamental lemma related to the error in  rank-1 approximations of tensors  of the form $\tens{X}+\tens{E}$, where $\tens{X}$ is a rank-1 tensor and $\tens{E}$ any other tensor, both with entries in some field $\mathbb{K}$.
\smallskip
\begin{lemma}\label{lemma:res}
Let $\tens{X}$ be a rank-1 tensor and $\phi$ the best rank-1 approximation operator. For any tensor $\tens{E}$,
$$
\|\tens{X}+\tens{E} - \phi(\tens{X}+\tens{E})\| \leq \sin(\gamma)\|\tens{E}\|,
$$ 
where  $\gamma$ denotes  the  angle between $\tens{E}$ and  $\tens{X}$.

\end{lemma}
\begin{proof}
Let $ \tens{P}_{\tens{X}}(\tens{X}+\tens{E})$ be the orthogonal projection of $\tens{X}+\tens{E}$ onto  $span(\tens{X})$. Because $\phi(\tens{X}+\tens{E})$ is a best rank-1 approximation of $\tens{X}+\tens{E}$, $ \tens{P}_{\tens{X}}(\tens{X}+\tens{E})$ cannot be a strictly better rank-1 approximation than $\phi(\tens{X}+\tens{E})$.  Thus, 
$$
\|\tens{X}+\tens{E} - \phi(\tens{X}+\tens{E})\| \leq \|\tens{X}+\tens{E}-\tens{P}_{\tens{X}}(\tens{X}+\tens{E})\|. 
$$ 
On the other hand,
 $\tens{X}+\tens{E} -  \tens{P}_{\tens{X}}(\tens{X}+\tens{E}) \perp \tens{X}$. Hence, we  have   $\|\tens{X}+\tens{E}-\tens{P}_{\tens{X}}(\tens{X}+\tens{E})\| = \sin(\gamma)\|\tens{E}\|$  by using basic trigonometry. 
 This concludes the proof. 
\end{proof}

\smallskip

The following results for the DCPD algorithm stems from the previous lemma. 
\smallskip

\begin{corollary}\label{cor:res}
The inequality $\|\tens{E}[r,l]\| \leq \sin(\gamma[r,l])\|\tens{E}[r-1,l]\|$ holds for any $1< r \leq R$. 
\end{corollary} 

\begin{proof}
By replacing $\tens{X}, \tens{E}$ and $\gamma$ in Lemma \ref{lemma:res} with $\tens{X}[r,l-1], \tens{E}[r-1,l]$ and $\gamma[r,l]$ respectively, the result follows directly.
\end{proof}
\smallskip

\begin{corollary}\label{cor:mds}
For any $l > 1$ and  $ c_{l} = \prod_{r=1}^{R}\sin(\gamma[r,l]),$ the inequality $\|\tens{E}[R,l]\| \leq c_{l} \|\tens{E}[R,l-1]\|$  holds.
\end{corollary} 

\begin{proof}
By applying $R-1$ times the result of Corollary \ref{cor:res}, we have
\[
\|\tens{E}[R,l]\| \leq \left(\sin(\gamma[R,l])\cdots \sin(\gamma[2,l])\right)\|\tens{E}[1,l]\|. 
\] From Lemma \ref{lemma:res}, we know that $$\|\tens{E}[1,l]\| \leq \sin(\gamma[1,l])\|\tens{E}[R,l-1]\|.$$ Thus, it follows that $\|\tens{E}[R,l]\| \leq c_{l}\|\tens{E}[R,l-1]\|$.
\end{proof}
\smallskip

Notice that the same result brought by Proposition (4.4) in  \cite{AlexIcassp2015} can be deduced from Corollary \ref{cor:mds} since $0 \leq c_{l} \leq 1$ for every iteration $l$, which implies the monotonic decrease of the sequence  $\{\|\tens{E}[R,l]\|\}_{l \in \mathbb{N}_{> 0}}$. 

\smallskip

\begin{lemma}\label{lemma:ortho}
If
$\|\tens{E}[R,l]\| = \|\tens{E}[R,l-1]\|$, then  $c_{l} = 1$.
\end{lemma}
\begin{proof}
Because $c_{l} \leq 1$ and  $\|\tens{E}[R,l]\| = \|\tens{E}[R,l-1]\|$, one concludes directly from Corollary \ref{cor:mds} that $c_{l}=1$.

\end{proof}

\smallskip
Lemma \ref{lemma:ortho} shows that the DCPD algorithm might not improve the estimation of the rank-1 components anymore for $l\geq l_{0}>1$.
And this may occur not only in the presence of noise.  
Actually, even for an almost orthogonal case $c_{l} \approx 1$,  $\|\tens{E}[R,l]\|$ may tend to a stationary non-zero value as $l$ increases. However, the DCPD algorithm converges to an exact decomposition if $c_{l} \leq C$, for all $l > 1$ for some constant $C<1$. This will be subsequently detailed by means of a geometric approach. 

\smallskip

Figure \ref{fig:cone} can also be seen as the representation of an $n$-sphere of dimension $n = I_1I_2\cdots I_N -1$ in $\mathbb{K}^{n+1}$ space. $\beta$ is half the white cone angle defined in $[0,\pi/2]$. The direction of the rank-1 tensor $\tens{X}[r,l-1]$ defines the axis of the white cone and varies with $r$ or $l$.  Under a condition on $\beta$, we can state an important proposition ensuring the convergence of the DCPD algorithm. 

\smallskip

\begin{proposition}\label{prop:beta}
Let $\tens{T}$ be a tensor such that $\rank{\tens{T}} \leq R$.  An exact decomposition is recovered by the DCPD algorithm if and only if there exists for every $(r,l)$ a half cone of angle $\beta$ (in white in Fig. \ref{fig:cone}), $0 \leq \beta < \pi/2$, such that $$\beta \geq \max_{l > 1}\min_{1\leq r\leq R}\gamma[r,l].$$  
\end{proposition}

\begin{proof}
$(\Leftarrow)$ For any iteration $l >1$, take $\gamma[r_{0},l] = \min_{1\leq r\leq R} \gamma[r,l].$ Notice that $c_{l} \leq \sin \left(\gamma[r_{0},l]\right)$ from Corollary \ref{cor:mds}. By hypothesis, $\sin \left(\gamma[r_{0},l]\right)\leq \sin(\beta)$, which implies that $ \|\tens{E}[R,l]\| \leq c_{l}\|\tens{E}[R,l-1]\| \leq \sin(\beta)\|\tens{E}[R,l-1]\|.$
 Because $\beta$ is an upper bound for $\gamma[r,l]$, $l >1$, we have
 $ \sin(\beta) \|\tens{E}[R,1]\| \geq \|\tens{E}[R,2]\|$ $\implies$ $\left(\sin(\beta)\right)^{l-1} \|\tens{E}[R,1]\|\geq \|\tens{E}[R,l]\|$.
Hence, when $l \rightarrow \infty$, $\|\tens{E}[R,l]\| \rightarrow 0$.

$(\Rightarrow)$ Let $l_{0}$ be some iteration such that $l_{0} > 1$. Without loss of generality, assume $\|\tens{E}[R,l]\| = 0$ for $l \geq l_{0}$ ($l_{0}$ can be arbitrarily large). Then $(\|\tens{E}[R,l]\|)_{l \in \mathbb{N}_{> 0}}$ is a strictly  monotonically decreasing sequence for $1 < l < l_{0}$, otherwise the algorithm would converge to a nonzero constant for some iteration smaller than $l_{0}$. Hence, for every $(r,l)$ we can choose $\beta$, $0 \leq \beta < \pi/2$ such that $\beta \geq \max \limits_{l > 1}\min \limits_{1\leq r\leq R}\gamma[r,l],$ and the proof is complete.
 \end{proof}

\smallskip

 As a conclusion, if for a given iteration $l$, all tensors $\tens{E}[r,l]$, $1\le r\le R$, fall within the gray volume depicted in Fig. \ref{fig:cone} (the complementary of the white cone), then the sequence $\tens{E}[r,l]$ does not tend to zero. Even if this gray volume can be made arbitrarily small, it is not of zero measure. So the best we can do is to prove an almost sure convergence of the DCPD algorithm to the exact decomposition  under some probabilistic conditions.

\smallskip

\begin{lemma} \label{lemma:cont}
If tensors $\tens{T}$ are distributed within $\set{T}^{(R)}$ according to  an absolutely  continuous probability measure, then $\|\tens{E}[r,l]\|$ are absolutely  continuous random variables.
\end{lemma}
\begin{proof}
Let $\vect{D} = [I_{1} \cdots I_{N}] $ be a specific size of $n$-order tensors, and let $\set{T}^{(R)}_{\vect{D}} = \{ \tens{T}  \in \mathbb{K}^{I_{1} \times \cdots \times I_{N}}:~  \tens{T} \subset \set{T}^{(R)} \}$. Because $\set{T}^{(R)} \supset \set{T}^{(R)}_{\vect{D}}$, any tensor $\tens{T}$ within $\set{T}^{(R)}_{\vect{D}}$ is also distributed according to  an absolutely  continuous probability measure.
Via the DCDP algorithm, each rank-1 component obtained in successive deflations is also in  $\set{T}^{(R)}_{\vect{D}}$. Hence, since the sum (subtraction) of continuous random variables does not affect the continuity, the residuals $\tens{E}[r,l]$  are also absolutely  continuous random variables.  Since  the norm is a $C^{0}$ function in finite dimension, $\|\tens{E}[r,l]\|$ is also absolutely continuous. 
\end{proof}
\smallskip

For the next developments, let $Z_{l} = \|\tens{E}[R, l ]\|$ and define the following probability  for some iteration $L>1$:
$$ 
F_{L}[\beta] =P\left (Z_{L} \leq \sin(\beta)Z_{L-1} \leq \ldots \leq \left(\sin(\beta)\right)^{L-1}Z_{1}\right). 
$$ 

$F_{L}[\beta]$ can be viewed as the probability that residuals fall within  at least one of the $R$ white cones in every iteration $l \leq L$. 

The following proposition ensures a reduction of  $Z_1$  by a factor smaller than  $\left(\sin(\beta)\right)^{L-1}$ after $L$ iterations with high probability, if a condition on the continuity of $F_{L}[\beta]$ is assumed.

\smallskip

\begin{proposition}\label{prop:contF}
Let $L$ be fixed. If  $\exists \beta_{0}: \beta_{0} \in [0, \pi/2)$ such that $F_{L}[\beta]$ is continuous on $[\beta_{0}, \pi/2]$, then $\forall \varepsilon:  \varepsilon \in (0,1],  \exists \beta \in [\beta_{0}, \pi/2)$ such that $F_{L}[\beta] > 1-\varepsilon$.
\end{proposition}
\smallskip

\begin{proof}
 Since $F_{L}[\pi/2]=1$ and  $F_{L}[\beta]$ is continuous on $[\beta_{0}, \pi/2]$, the proof follows directly from the intermediate value theorem. \end{proof}
\smallskip

Although $Z_{l}, 1 \leq l \leq L,$ are absolutely continuous random variables and $g_{m}(\beta) = \sin^{m}(\beta)$ are continuous functions for all $m \geq 0$, the continuity of  $F_{L}[\beta]$ in $\beta$ is not  guaranteed   due to the dependence of the random variables $Z_{1}, \ldots, Z_{L}$ ($Z_{l}$ depends on $Z_{l-1}$). For example, for $L = 2$ and $Z_{1} = 2Z_{2}$ with probability $1$, it is easy to check that $F_{2}[\beta]$ is not continuous at $\beta = \pi/2$. Indeed, $\lim_{\beta \rightarrow \pi/2^{-}} F_{2}[\beta] = 0$ whereas $F_{2}[\pi/2] = 1$.

The following conjecture claims that there exists absolutely continuous distributions of tensors in $\set{T}^{(R)}$ such that the probability $F_{l}[\beta(l)]$ tends to $1$ for some function $\beta(l)$ as  $l \rightarrow \infty$,  and at the same time the norm of residuals tends to $0$, which is suitable for the convergence of the DCPD algorithm to an exact CP decomposition.
  
\smallskip

\begin{conjecture}\label{conj:prob}
There exists at least one absolutely continuous probability measure $\mu$ for tensors $\tens{T}$ within $\set{T}^{(R)}$ for which the following holds: 
\begin{itemize}
\item[(i).]~$\forall \varepsilon: \varepsilon \in (0,1]$, $\forall l: l > 1$,  $\exists \beta: \beta \in [0, \pi/2)$,  such that $F_{l}[\beta] > 1-\varepsilon$.
\item[(ii).]~$\forall l: l > 1$,  $\exists \beta: \beta \in [0, \pi/2)$, such that  $\left(\sin(\beta)\right)^{l-1}$ is a strictly monotonically decreasing sequence converging to $0$.
\end{itemize}
\end{conjecture}
  \smallskip

Subsequent computer simulations support the existence of  a uniform probability measure $\mu$ for the entries of tensors within $\set{T}^{(R)}$, such that $F_{\beta} [L] \approx 1$ for  large values of $L$,  and $\|\tens{E}[R,L]\| \approx 0$.  This reinforces our conjecture.

\section{Computer Results}\label{sec:sim_res}

\subsection{Comparison between THOSVD and SeROAP}\label{sec:TS}
\begin{figure}[h!] 
   \centering
	   \psfrag{T1}[cc][cc]{ \small $3 \times 4 \times 5$ tensors} 
	     \psfrag{T2}[cc][cc]{ \small $3 \times 4 \times 20$ tensors} 
	       \psfrag{T3}[cc][cc]{ \small $3 \times 20 \times 20$ tensors} 
	         \psfrag{T4}[cc][cc]{ \small $20 \times 20 \times 20$ tensors} 
	          \psfrag{Y1}[cc][cc]{  \footnotesize $\Delta\phi$} 
	           \psfrag{X1}[cc][cc]{ } 
	            \psfrag{X2}[cc][cc]{ } 
	    \psfrag{Y2}[cc][cc]{  \footnotesize  $\Delta\phi$} 
	     \psfrag{Y3}[cc][cc]{  \footnotesize$\Delta\phi$} 
	      \psfrag{Y4}[cc][cc]{   \footnotesize $\Delta\phi$} 
	    ~\hspace{-1.4em}\includegraphics[width=0.45\textwidth]{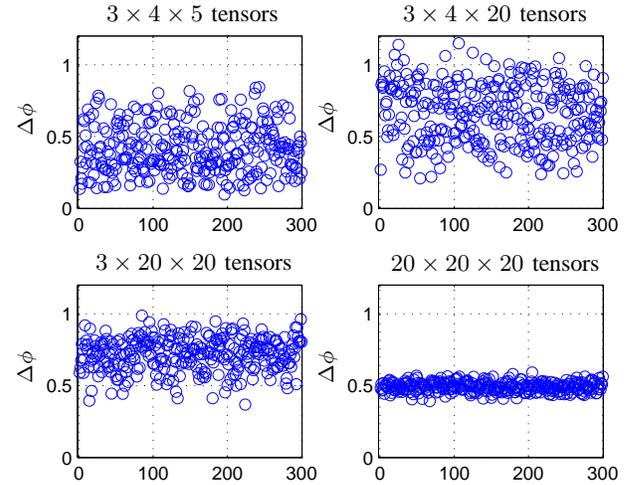}
  \vspace{-2.3ex}
     \caption{THOSVD and SeROAP comparison}\label{fig:TxS}
\end{figure}

 In this section,  we compare the performance of  rank-1 approximation methods SeROAP and THOSVD for different three-way tensor scenarios. For each case, $300$ complex tensors whose real and imaginary parts  are uniformly distributed in $[-1,1]$ were generated.  Figure \ref{fig:TxS} presents the difference between the Frobenius norms  of the residuals computed as $$\Delta\phi = \|\tens{T}-\phi^{\textrm{TH}}(\tens{T})\| -  \|\tens{T} - \phi^{\textrm{Se}}(\tens{T})\|.$$

We note that $\Delta\phi > 0$ in all scenarios, as predicted by Proposition  \ref{prop:THSe}.

\subsection{Performance of rank-1 approximations}\label{sec:rank_one}

The tables below  compare different rank-1 approximation methods with respect to the best rank-1 approximation, which was obtained from the algebraic geometric moment method described in \cite{lasserre2001global}. Because the latter  is  infeasible to compute for high dimensions, we have focused on $2\times 2\times 2$ and $3\times 3\times 3$ real tensors.

The two iterative methods, namely ALS and CE, are initialized by the result obtained from SeROAP. A sample of $200$ real tensors uniformly distributed in $[-1,1]$ were generated for each of both scenarios. For comparison, we consider the MSE metric given by

\[
\text{MSE} = \frac{1}{200} \sum \limits_{n=1}^{200} \left (\Delta\phi_{m}^{(n)} \right)^{2},
\]
where $\Delta\phi_{m}^{(n)} = \|\tens{T}^{(n)}-\phi_{m}(\tens{T}^{(n)})\| -  \|\tens{T}^{(n)} - \phi^{\star}(\tens{T}^{(n)})\|$. $\phi_{m}(\tens{T}^{(n)})$ and $\phi^{\star}(\tens{T}^{(n)})$ are the rank-1 approximation for algorithm $m$ and the best rank-1 approximation of $\tens{T}^{(n)}$, respectively.

\begin{center} 
\begin{tabular}{|c||c|c|} 
\cline{2-3}\multicolumn{1}{c|}{} & \multicolumn{2}{c||}{$2\times 2 \times 2$ tensors} \\ \hline 
Algorithm & MSE & mean iteration \\ \hline \hline
THOSVD  & 0.02299 &N/A   \\ \hline
SeROAP &4.36155e-4  &N/A    \\ \hline
CE &1.13056e-18 &6.115   \\ \hline 
ALS & 1.03406e-13 &6.135  \\ \hline 
\end{tabular} 
\end{center}

\begin{center} 
\begin{tabular}{|c||c|c|} 
\cline{2-3}\multicolumn{1}{c|}{} & \multicolumn{2}{c||}{$3\times 3 \times 3$ tensors} \\ \hline 
Algorithm & MSE & mean iteration  \\ \hline \hline
THOSVD  & 0.08386 & N/A   \\ \hline
SeROAP  &0.00172 & N/A   \\ \hline
CE &1.70946e-4  &10.990 \\ \hline 
ALS  &1.70949e-4  &10.435 \\ \hline 
\end{tabular} 
\end{center}
The results show that SeROAP is a better rank-1 approximation than THOSVD as expected. For $2\times 2\times 2$ tensors, CE attains the best rank-1 approximation. In both scenarios, ALS and CE converge approximately in the same number of iterations.

\subsection{Percentage of successful decompositions}\label{sec:result1}

\begin{figure}[ht!] 
 \small
   \centering
   \psfrag{TTT}[cc][cc]{} 
   \psfrag{A}[cc][cc]{\footnotesize \text{3x3x3}} 
    \psfrag{B}[cc][cc]{\footnotesize 4x4x4} 
     \psfrag{C}[cc][cc]{\footnotesize 5x5x5} 
      \psfrag{D}[cc][cc]{\footnotesize 6x6x6} 
       \psfrag{E}[cc][cc]{\footnotesize 7x7x7} 
        \psfrag{F}[cc][cc]{\footnotesize 8x8x8} 
  ~\hspace{-1.4em}\includegraphics[width=0.4\textwidth]{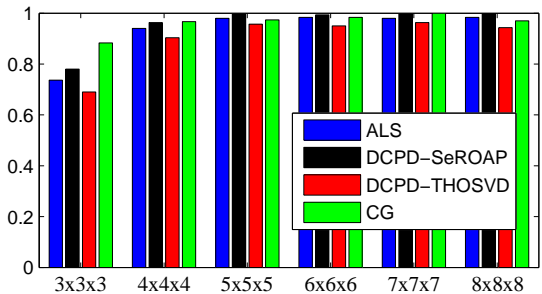}
  \vspace{-2.0ex}
     \caption{Percentage of successful decomposition for rank-$3$ tensors.}\label{fig:percent}
\end{figure}

Figure \ref{fig:percent} presents the percentage of successful decompositions of rank-$3$ tensors for the algorithms ALS, CG with ELS, and DCPD.  The ALS and the CG algorithms were randomly initialized. We have simulated DCPD with the algebraic methods THOSVD and SeROAP. Noise is not considered in this case so that the performance is evaluated for  the computation of  an  exact decomposition of $300$ tensors. We consider that a decomposition is succeeded if the residual $\|\tens{E}\| \leq 10^{-6}$.

We note that the DCPD algorithm combined with SeROAP always presents a better performance than the  standard  ALS algorithm. Moreover, for higher dimensions, the percentage of successful decompositions is almost  $100\%$ for DCPD-SeROAP, which is a remarkable result, bearing in mind that the objective is multimodal.

\subsection{ Convergence rate }

\begin{figure} 
\centerline{ \hspace{-2ex}
\includegraphics[width=0.5\textwidth]{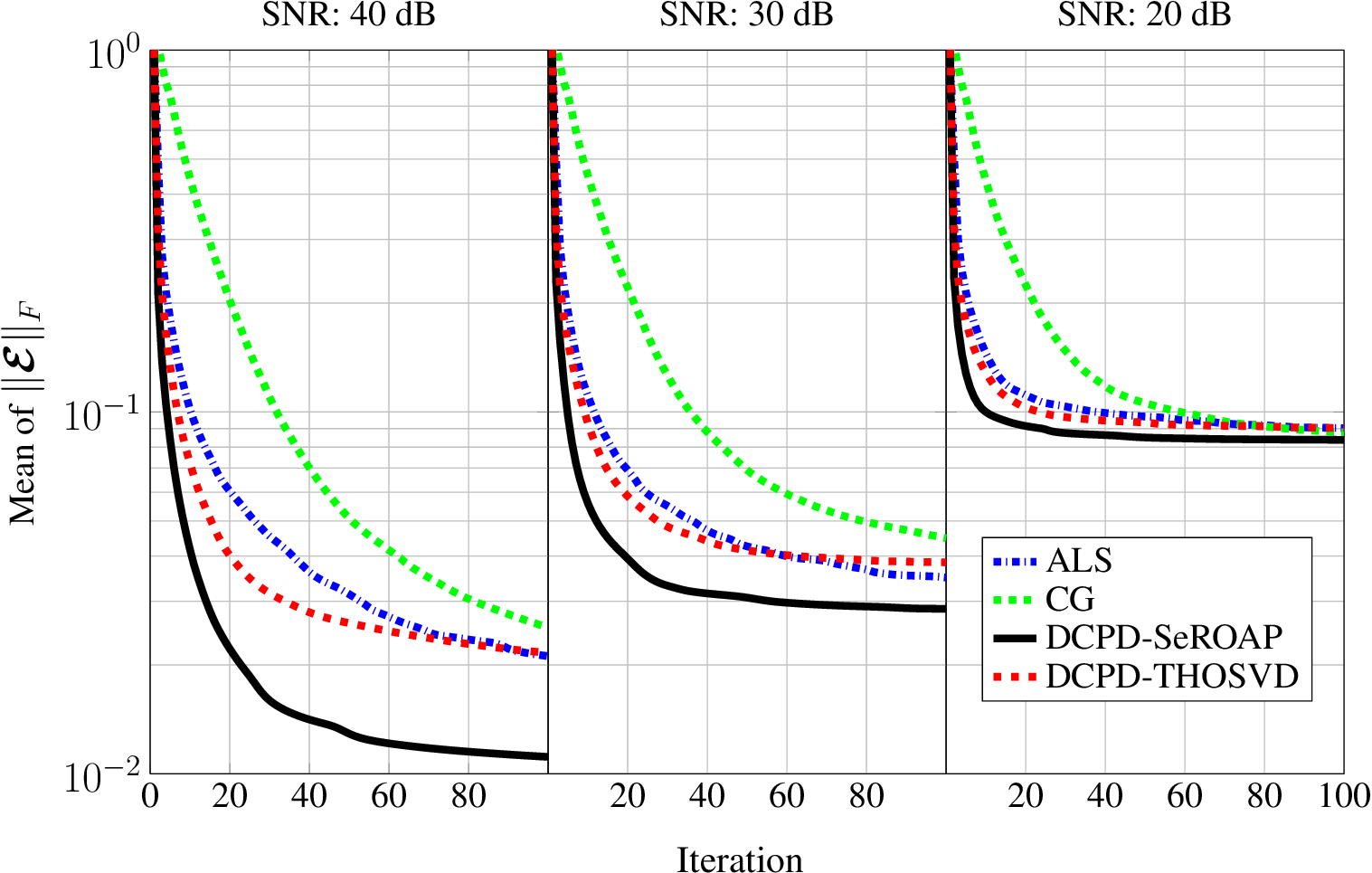}}
 \vspace{-2.0ex}
\caption{Mean of $\|\tens{E}[R,l]\|$ per iteration for different values of SNR.}\label{fig:ErrorItr}
 \end{figure}

Figure \ref{fig:ErrorItr} presents the performance of the algorithms in terms of the average of  $\|\tens{E}[R,l]\|$  per iteration for different values of the signal-to-noise (SNR) ratio for $5\times5\times5$- rank-$3$ tensors. Again,  ALS and  CG  algorithms were randomly initialized. Additive  Gaussian   noise is considered in our simulations.

We note  in this figure that DCPD-SeROAP converges more quickly  than the other algorithms. For an SNR of $40$ dB, DCPD-SeROAP attains   $\|\tens{E}[R,l]\| \approx 0.01$  in approximately $100$ iterations while, for the other algorithms,  $\|\tens{E}[R,l]\| > 0.02$  for the same number of iterations. Similar results are observed for other SNRs. The figure also shows that performances become similar  when the SNR is decreased.

\subsection{Residual vs rank }
Now, we compare the algorithms for two SNRs by varying the rank of $8 \times 8\times 8$ tensors. 
\begin{figure}[ht!]
	 \large
	   \centering
	   \psfrag{T1}[cc][cc]{\footnotesize $8\times8\times8:$ $\text{SNR} = 40 \hspace{1mm} \text{dB}$ } 
	    \psfrag{T2}[cc][cc]{\footnotesize $8\times8\times8:$ $\text{SNR} = 30 \hspace{1mm}  \text{dB}$ } 
	   \psfrag{A1}[cc][cc]{\footnotesize rank-$3$} 
	    \psfrag{B1}[cc][cc]{\footnotesize rank-$3$} 
	   \psfrag{A2}[cc][cc]{\footnotesize rank-$4$} 
	    \psfrag{B2}[cc][cc]{\footnotesize rank-$4$} 
	    \psfrag{A3}[cc][cc]{\footnotesize rank-$5$} 
	    \psfrag{B3}[cc][cc]{\footnotesize rank-$5$} 
	    \psfrag{A4}[cc][cc]{\footnotesize rank-$6$} 
	    \psfrag{B4}[cc][cc]{\footnotesize rank-$6$} 
	    \psfrag{A5}[cc][cc]{\footnotesize rank-$7$} 
	    \psfrag{B5}[cc][cc]{\footnotesize rank-$7$} 
	   \psfrag{Y1}[cc][cc]{\footnotesize Mean of $\|\tens{E}\|$} 
	   \psfrag{Y2}[cc][cc]{\footnotesize Mean of $\|\tens{E}\|$} 
	   ~\hspace{-1.4em}\includegraphics[width=0.5\textwidth]{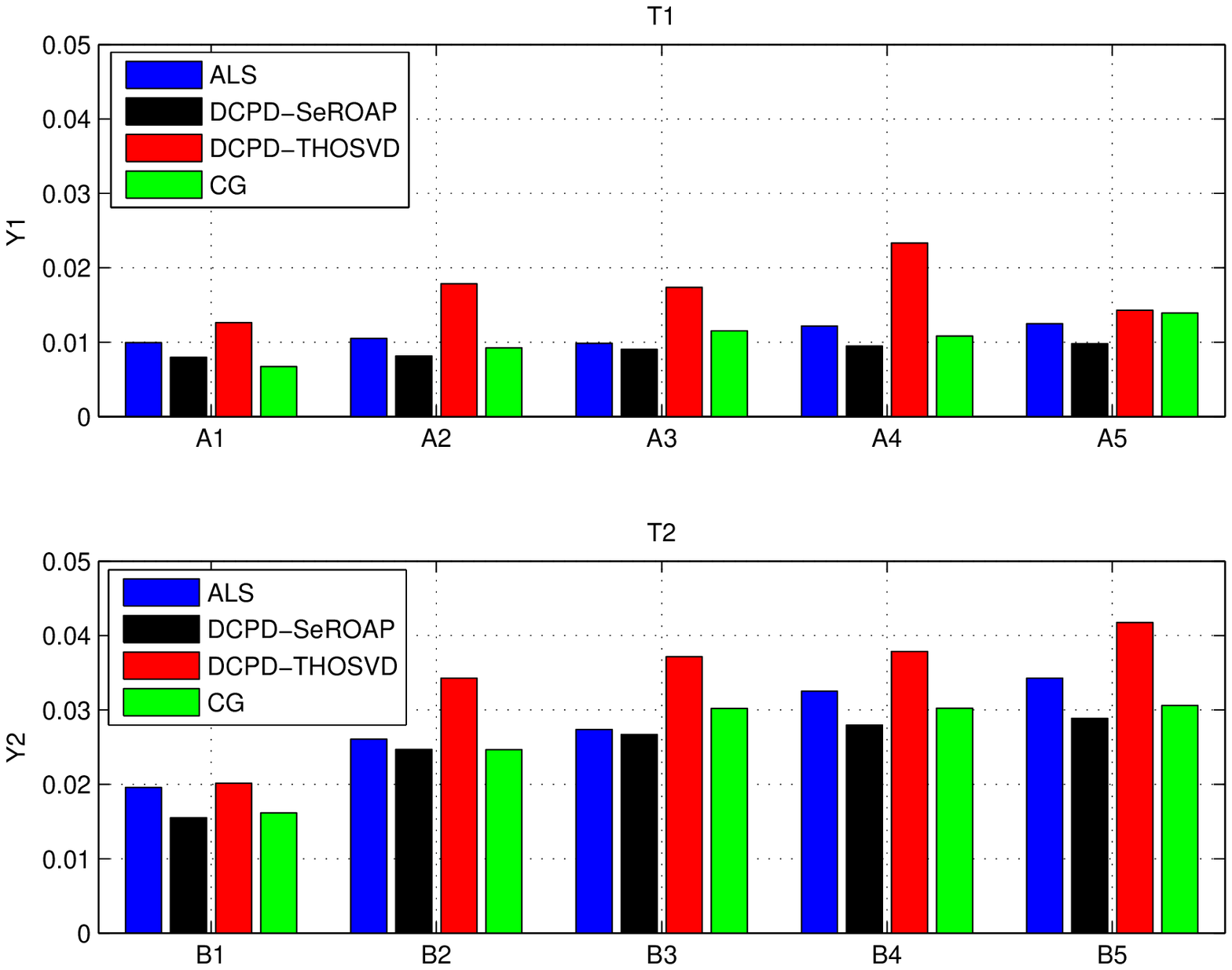}
	  \vspace{-2.0ex}
	     \caption{Mean of $\|\tens{E}\|$ under rank variation.}\label{fig:rankvar}
	\end{figure}
Again, we note the better performance of DCPD-SeROAP over the competing algorithms.  Figure \ref{fig:rankvar} also shows that the combination of DCPD and THOSVD yields the worst results. This is expected because the rank-1 approximation obtained by THOSVD is not good enough, so that  DCPD does not come at a small residual.

\section{Conclusion}\label{sec:concl}
In this paper, we presented some CP tensor decompositions algorithms and provided an analysis of their computational complexities. Our contributions included:  (i) a new algebraic rank-1 method, namely SeROAP, performing better than THOSVD for three-way tensors; (ii) an iterative rank-1 approximation  algorithm, namely CE,  that  refines  any rank-1 approximation method, such as SeROAP and THOSVD,  which converges  in very few iterations; and (iii) an  analysis of the convergence of the DCPD algorithm under a geometric point of view.
Several computer experiments have confirmed the theoretical results.

\appendices
\section{Computation of rank-1 approximation using SeROAP} \label{app:SeROAP}

We present an example of how SeROAP algorithm works for computing a rank-1 approximation of a given tensor. 

Let $\tens{T}$ be a $2\times2\times2 \times 2$ complex tensor whose mode-$1$ unfolding is given by 

$$\matr{T}_{(1)} = \begin{bmatrix}
1 & -1 & 0 & 1 & 3 & i & 0 & 1\\
0 & 1 & -i & 1 & 1 & 0 & 2 & -2i
\end{bmatrix}.$$
In the first \textit{for} loop of SeROAP algorithm, $\vect{v}_{1}$ is the dominate right singular vector of $\matr{V} = \matr{V}_{0} = \matr{T}_{(1)}$ given by

$$\vect{v}_{1} = \begin{bmatrix}
 -0.1717 - 0.0914i\\
   0.0245 + 0.1060i\\
  -0.0146 - 0.1472i\\
  -0.3189 - 0.0768i\\
  -0.6624 - 0.2596i\\
  -0.0914 + 0.1717i\\
  -0.2944 + 0.0292i\\
  -0.2010 - 0.3858i
\end{bmatrix}.
$$  By reshaping $\vect{v}_{1}$ in a $2\times4$ matrix $\matr{V}_{1}$ we have


\begin{align*}
 \matr{V}_{1}  =   & \left [ \begin{matrix}
-0.1717 - 0.0914i&   -0.0146 - 0.1472i\\
0.0245 + 0.1060i &   -0.3189 - 0.0768i
\end{matrix} \right. & \\
  & \qquad \left. \begin{matrix}
 -0.6624 - 0.2596i & -0.2944 + 0.0292i\\
  -0.0914 + 0.1717i &   -0.2010 - 0.3858i
\end{matrix} \right ].
\end{align*}

In next iteration ($n=2$), we compute the dominate right singular vector of $\matr{V} = \matr{V}_{1}$. Hence, 
$$\vect{v}_{2} =  \begin{bmatrix}
   -0.1654 + 0.1657i\\
   0.0611 + 0.2758i\\
  -0.6190 + 0.6261i\\
  -0.2033 + 0.2210i
\end{bmatrix}, 
$$ and $\matr{V} = \matr{V}_{2} = Unvec(\vect{v}_{2})$ is updated 
$$
 \matr{V}  = \begin{bmatrix}
-0.1654 + 0.1657i &   -0.6190 + 0.6261i \\
 0.0611 + 0.2758i  &  -0.2033 + 0.2210i
 \end{bmatrix}.
$$

The next step is to compute the vector $\vect{w} = \vect{v}^{*} \kron \vect{u}$, where $\vect{u}$ and $\vect{v}$ are the first left and right singular vectors, respectively. Thus,
$$\vect{u} =  \begin{bmatrix}
  0.6106 - 0.7024i\\
   0.1758 - 0.3208i
\end{bmatrix}, 
\vect{v} =  \begin{bmatrix}
    -0.3027 - 0.0545i \\
  -0.9481 + 0.0809i
   \end{bmatrix},
  $$ 
  $$
   \vect{w} =  \begin{bmatrix}
  -0.1466 + 0.2459i \\
  -0.0357 + 0.1067i\\
  -0.6357 + 0.6166i\\
  -0.1926 + 0.2899i
  \end{bmatrix}.
$$

Notice that $\vect{w}$ can be viewed as a vectorization of a rank-$1$ matrix. In the end of the first iteration of the second \textit{for} loop, $\vect{w}$ is updated so that it becomes a vectorization of a rank-$1$ three-way tensor. This is achieved by performing the projection of each row of  matrix $\matr{V}_{1}$ onto $\vect{w}$, by means of $\matr{V}_{1}\vect{w}$, following by the multiplication with $\vect{w}^{H}$. Hence
\begin{align*}
 \matr{X}_{(2)}  =   & \left [ \begin{matrix}
 -0.1900 - 0.1178i & -0.0631 - 0.0611i\\
-0.0604 + 0.0051i & -0.0236 - 0.0031i
\end{matrix} \right. & \\
  & \qquad \left. \begin{matrix}
-0.6624 - 0.1989i &  -0.2377 - 0.1317i\\
-0.1762 + 0.0638i & -0.0730 + 0.0098i
\end{matrix} \right ],
\end{align*} and 


$$\vect{w} =  \begin{bmatrix}
   -0.1900 - 0.1178i\\
  -0.0604 + 0.0051i\\
  -0.0631 - 0.0611i\\
  -0.0236 - 0.0031i\\
  -0.6624 - 0.1989i\\
  -0.1762 + 0.0638i\\
  -0.2377 - 0.1317i\\
  -0.0730 + 0.0098i
\end{bmatrix}
$$

In the next iteration, we perform $\matr{X}_{1} = (\matr{V}_{0})\vect{w}^{H}$, which is an unfolding of a  $4$-order rank-$1$  tensor.
Actually, in every iteration of the second \textit{for} loop, $\matr{X}_{n}$ is updated to  an  unfolding of a rank-$1$ tensor of order $N-n+1$. Hence,

\begin{align*}
 \matr{X}_{(1)}  =   & \left [ \begin{matrix}
0.5373 - 0.0992i  &  0.1329 + 0.0653i\\
  0.2696 - 0.1010i &  0.0750 + 0.0216i
\end{matrix} \right. & \\
  & \qquad \left. \begin{matrix}
 0.1981 - 0.0830i &  0.0565 + 0.0140i\\
 0.0951 - 0.0613i &  0.0306 + 0.0020i
\end{matrix} \right. & \\  
  & \qquad \quad \left. \begin{matrix}
1.6851 + 0.1360i &  0.3446 + 0.3020i \\
0.8868 - 0.0849i  & 0.2066 + 0.1249i
\end{matrix} \right. & \\ 
  & \qquad \qquad \left. \begin{matrix}
0.6585 - 0.0886i &  0.1576 + 0.0872i \\
 0.3335 - 0.1067i &  0.0898 + 0.0307i
\end{matrix} \right ],
\end{align*}

Indeed, $\matr{X}_{(1)}$ is the mode-$1$ unfolding of the rank-$1$ approximation of $\tens{T}$ computed by SeROAP algorithm.

\bibliographystyle{IEEEbib}
\bibliography{biblioAlg}

\end{document}